\documentclass[a4paper,11pt,twoside,reqno]{amsart}

\usepackage[utf8]{inputenc}
\usepackage[plainpages=false,pdfpagelabels=true]{hyperref}
\usepackage{amsmath}
\usepackage{amssymb,amsthm}
\usepackage[margin=1in]{geometry}
\usepackage{slashed}

\newtheorem{Satz}{Theorem}[section]
\newtheorem{Prop}[Satz]{Proposition}
\newtheorem{Lem}[Satz]{Lemma}

\newtheorem{Cor}[Satz]{Corollary}
\theoremstyle{definition}
\newtheorem{Dfn}[Satz]{Definition}
\newtheorem{Bem}[Satz]{Remark}

\newcommand{\vol}{{\operatorname{vol}}}

\newcommand{\s}{{\mathbb{S}}}
\allowdisplaybreaks[1]

\renewcommand{\epsilon}{\varepsilon}

\newcommand{\R}{\ensuremath{\mathbb{R}}}

\numberwithin{equation}{section}

\title{On p-biharmonic curves}
\author{Volker Branding}
\date{\today}
\address{University of Vienna, Faculty of Mathematics\\
Oskar-Morgenstern-Platz 1, 1090 Vienna, Austria\\}
\email{volker.branding@univie.ac.at}
\subjclass[2010]{58E20; 53C43; 31B30; 58E10}
\keywords{p-biharmonic curves; p-elastic curves; stability}
\thanks{The author gratefully acknowledges the support of the Austrian Science Fund (FWF) through the START-project (Y963) of Michael Eichmair and 
the project "Geometric Analysis of Biwave Maps" (P34853).
}
\begin{document}

\begin{abstract}
In this article we study \(p\)-biharmonic curves as a natural generalization
of biharmonic curves. In contrast to biharmonic curves \(p\)-biharmonic
curves do not need to have constant geodesic curvature if \(p=\frac{1}{2}\)
in which case their equation reduces to the one of \(\frac{1}{2}\)-elastic curves.
We will classify \(\frac{1}{2}\)-biharmonic curves on closed surfaces
and three-dimensional space forms making use of the results obtained
for \(\frac{1}{2}\)-elastic curves from the literature.
By making a connection to magnetic geodesic we are 
able to prove the existence of
\(\frac{1}{2}\)-biharmonic curves on closed surfaces.
In addition, we will discuss the stability of \(p\)-biharmonic curves
with respect to normal variations.
Our analysis highlights some interesting relations between \(p\)-biharmonic
and \(p\)-elastic curves.
\end{abstract} 

\maketitle

\section{Introduction and Results}
Finding interesting curves on a Riemannian manifold is one of the central
topics in modern differential geometry. The most prominent examples of
such curves are without doubt \emph{geodesics} as they are the curves with minimal
distance between two given points in a Riemannian manifold.
Moreover, the existence of geodesics is guaranteed by a number of famous
results in differential geometry, such as the Theorem of \emph{Hopf-Rinow}.

In order to obtain the equation for geodesics it is useful
to employ a variational approach as this gives a lot of 
additional mathematical structure.
To this end, we consider a curve \(\gamma\colon I\to M\),
where \(I\subset\R\) represents an interval, \((M,g)\) a Riemannian manifold
and by \(s\) we denote the parameter of the curve \(\gamma\).
Moreover, we use the notation \(\gamma'=\frac{d\gamma}{ds}\).
Then, we define the energy of the curve \(\gamma\) by
\begin{align}
\label{energy-curve}
E_1(\gamma)=E(\gamma):=\frac{1}{2}\int_I|\gamma'|^2ds.
\end{align}

The critical points of \eqref{energy-curve} are precisely \emph{geodesics}
and are characterized by the equation
\begin{align*}
\tau(\gamma):=\nabla_{\gamma'}\gamma'=0,
\end{align*}
which is a second order non-linear ordinary differential equation for \(\gamma\)
and the quantity \(\tau(\gamma)\) is usually referred to as \emph{tension field}.
In the case of a higher-dimensional domain 
\eqref{energy-curve} is replaced by the energy of a map
whose critical points are \emph{harmonic maps}.

Another interesting class of curves can be obtained by 
extremizing the \emph{bienergy} of a curve \(\gamma\)
which is given by 
\begin{align*}
E_{2}(\gamma):=\frac{1}{2}\int_I|\tau(\gamma)|^2ds.
\end{align*}
The critical points of this functional are called \emph{biharmonic curves}
and are characterized by a non-linear ordinary differential equation
of fourth order. 
However, as biharmonic curves necessarily have constant geodesic curvature
they are often too rigid in order to allow for applications
in the natural sciences. In order to circumvent this drawback
there exist several notions similar to biharmonic curves 
which are not that restricted.

In this article we will introduce one such further generalization of
biharmonic curves. The starting point of our analysis
is the energy for \emph{p-biharmonic curves} which is defined as follows
\begin{align}
\label{energy-p-biharmonic}
E_{2,p}(\gamma)=\frac{1}{p}\int_I|\tau(\gamma)|^pds,\qquad p>0.
\end{align}
This definition is motivated from applications of such curves
in elasticity where it is often favorable to choose an exponent 
different from \(2\). Of course, for \(p=2\) the energy functional
\eqref{energy-p-biharmonic} reduces to the well-studied energy for
biharmonic curves.

Another possible variant 
are the so-called \emph{p-elastic curves}
which are critical points of the \emph{p-elastic energy} given by
\begin{align}
\label{p-elastic-energy}
E^{elastic}_p(\gamma)=\frac{1}{p}\int_Ik^pds,\qquad p>0.
\end{align}
Here, \(k\) represents the geodesic curvature of the curve \(\gamma\).

Although both \eqref{energy-p-biharmonic}
and \eqref{p-elastic-energy} look very similar at first glance we can expect them to be different in general: If we could consider \eqref{energy-p-biharmonic}
and \eqref{p-elastic-energy} in the case of a two-dimensional domain with \(p=2\)
then \eqref{energy-p-biharmonic} would become the bienergy while
\eqref{p-elastic-energy} would turn into the Willmore energy. As these two variational
problems are very different in nature one should also expect substantial differences
between \eqref{energy-p-biharmonic} and \eqref{p-elastic-energy}.

From an analytic point of view both \eqref{energy-p-biharmonic}
and \eqref{p-elastic-energy} are favorable to investigate if \(p\geq 2\)
as the resulting Euler-Lagrange equations might become degenerate
if \(p<2\). One could also allow for the possibility of \(p\) being negative.
However, as we always would like to draw a comparison between critical points of
\eqref{energy-p-biharmonic}, \eqref{p-elastic-energy} and geodesics we restrict to positive values of \(p\) as for \(p<0\) we need to be careful with \eqref{energy-p-biharmonic}, \eqref{p-elastic-energy} should the curve be a geodesic.

Let us give an (non-exhaustive) overview on the mathematical results
on \(p\)-elastic curves.
One of the most influential works on this subject for the case \(p=2\)
goes back to Langer and Singer 
\cite{MR772124}
who classified elastic curves on two-dimensional manifolds.
Their work was later extended by Watanabe 
\cite{MR3229087}
who was able to determine the geodesic curvature of \(p\)-elastic
curves for \(p\geq 2\) making use of generalized elliptic integrals.
Recently, Shioji and Watanabe showed some interesting
phenomena for \(p\)-elastic curves on \(\s^2\) in
\cite{MR4184824}.
Concerning the gradient flow of \(p\)-elastic curves for \(p\geq 2\)
we refer to the recent articles
\cite{MR4418802,MR4150936,MR4318851}. The elastic flow (\(p=2)\)
of curves on the sphere was investigated in \cite{MR3781340}.
For the current status of research regarding the flow of elastic
curves we refer to the recent survey \cite{MR4277362}.

One of the main observations of this article is
that in the case of \(p=\frac{1}{2}\) the first variations
of \eqref{energy-p-biharmonic} and \eqref{p-elastic-energy}
lead to the same equation. While for \(p\neq\frac{1}{2}\)
the equations for \(p\)-biharmonic curves always require
that we have a curve of constant curvature we will see that
for \(p=\frac{1}{2}\) we do not encounter this restriction.

The idea of studying \(p\)-elastic curves goes back to Bernoulli,
and later Blaschke initiated the study of \(\frac{1}{2}\)-elastic curves in \(\R^3\).
Recently, \(\frac{1}{2}\)-elastic curves have been investigated
in the case that the ambient space is a sphere or hyperbolic 
space in a series of articles by Arroyo, Garay and Menc\'ia
\cite{MR2007599,MR2078687} and by Arroyo, Garay and Barros \cite{MR2031956}.
Another variant of \eqref{p-elastic-energy} that received growing
attention is
\begin{align}
\label{energy-blaschke}
E_{\mu}(\gamma)=\int_I\sqrt{k+\mu}ds,\qquad \mu\in\R,
\end{align}
Recently the functional \eqref{energy-blaschke} has been investigated 
in great generality by Arroyo, Garay and P\'ampano \cite{MR3774309,MR4028018}.
At various places our analysis is close to the
results obtained in \cite{MR2031956,MR2007599,MR2078687,MR3774309,MR4028018}
in the context of critical points of 
\eqref{energy-blaschke}.
A general discussion on curvature energies of curves can be found in 
\cite{MR2076749}.

In the case of a higher-dimensional domain \eqref{energy-p-biharmonic}
becomes the energy for \(p\)-biharmonic maps. A number of geometric
classification results for the latter could be achieved in 
\cite{MR3567234,MR4169439,MR3393120}.
Some stability results on F-biharmonic maps, which are
a further generalization of \(p\)-biharmonic maps, 
have been obtained in
\cite{MR3773492}. Recently, \(p\)-biharmonic submanifolds
were investigated in \cite{MR4594934}.
For the stability of biharmonic maps we recommend the reader
to consult the recent articles \cite{MR4216418,MR4110268,MR4386842} and the references therein. A stability analysis for harmonic 
self-maps of cohomogeneity one manifolds, which employs similar ideas
as the ones used in this article,
was carried out in \cite{MR4477489}.
Concerning the current status of research on biharmonic curves we refer to \cite[Section 4]{MR3098705}, a general introduction to the field 
of biharmonic maps in Riemannian geometry is provided by the recent book \cite{MR4265170}.

Let us also mention the following results closely connected to this article.
The current status of research on higher order variational problems 
can be found in \cite{MR4106647}. A family of curves that interpolates between geodesics and biharmonic curves
has been investigated in 
\cite[Section 3]{MR4058514}.

Throughout this article we will use the following notation:
We use \(s\) to represent the parameter of the curve \(\gamma\)
and denote differentiation of the curve \(\gamma\) with respect
to the curve parameter by \(\gamma'\). By \(I\subset\R\) we denote
an interval and \((M,g)\) represents a Riemannian manifold of dimension \(\dim M=n\).
For the Riemannian curvature tensor we use the sign convention
\begin{align*}
R(X,Y)Z=\nabla_X\nabla_YZ-\nabla_Y\nabla_XZ-\nabla_{[X,Y]}Z,
\end{align*}
where \(X,Y,Z\) are vector fields.

This article is organized as follows:
In Section 2 we study the first and second variation of the energy functional
for \(p\)-biharmonic curves and discuss the relation between \(p\)-biharmonic
and \(p\)-elastic curves.
Afterwards, in Section 3, we turn to the analysis
of \(p\)-biharmonic curves on surfaces. 
In particular, we study their 
stability with respect to normal variations.
Moreover, employing a number of results on magnetic geodesics from 
the literature we are able to obtain some existence results
for \(p\)-biharmonic curves.
Section 4 is devoted to the analysis of \(p\)-biharmonic curves on
three-dimensional space forms also including an analysis of their stability.
% Finally, in Section 5, we derive the equation for \(p\)-biharmonic curves on the Euclidean sphere of arbitrary dimension.

\section{Variational formulas}
In this section we derive the first and second variation of \eqref{energy-p-biharmonic} and discuss the relation between \(p\)-biharmonic and \(p\)-elastic
curves.
\subsection{The first variation formula}
We start by deriving the Euler-Lagrange equation of \eqref{energy-p-biharmonic}.
To this end, we fix a small number \(\epsilon>0\) and let \(\gamma_t\colon (-\epsilon,\epsilon)\times I\to M\) be a variation
of the curve \(\gamma\) satisfying \(\frac{\nabla\gamma_t}{\partial t}\big|_{t=0}=\eta\).
For the moment we assume that \(\eta\in\Gamma(\gamma^\ast TM)\).
Here, \(I\subset\R\) is a closed interval.

\begin{Lem}
Let \(\gamma_t\) be a variation  of \(\gamma\) as described above. 
Then the following formula holds
\begin{align}
\label{first-variation-p-bienery}
\frac{d}{dt}\big|_{t=0}\frac{1}{p}\int_I|\nabla_{\gamma_t'}\gamma_t'|^pds=
\int_I\bigg(\langle \eta,
\nabla_{\gamma'}\nabla_{\gamma'}\big(|\nabla_{\gamma'}\gamma'|^{p-2}\nabla_{\gamma'}\gamma'\big)
 -|\nabla_{\gamma'}\gamma'|^{p-2}R^M(\gamma',\nabla_{\gamma'}\gamma')\gamma'
\rangle \bigg)ds.
\end{align}
\end{Lem}
\begin{proof}
It is straightforward to show that
\begin{align*}
\frac{\nabla}{\partial t}\big|_{t=0}\nabla_{\gamma_t'}\gamma_t'=\nabla_{\gamma'}\nabla_{\gamma'}\eta
+R^M(\eta,\gamma')\gamma',
\end{align*}
where we first interchanged covariant derivatives and then 
employed the torsion-freeness of the Levi-Civita connection.

Then, we get
\begin{align*}
\frac{d}{dt}\big|_{t=0}\frac{1}{p}\int_I|\nabla_{\gamma_t'}\gamma_t'|^pds
=&\int_I|\nabla_{\gamma'}\gamma'|^{p-2}\big(
\langle\nabla_{\gamma'}\nabla_{\gamma'}\eta,\nabla_{\gamma'}\gamma'\rangle
+\langle R^M(\eta,\gamma')\gamma',\nabla_{\gamma'}\gamma'\rangle
\big)ds.
\end{align*}
Finally, using the symmetries of the Riemannian curvature tensor
and integration by parts completes the proof.
\end{proof}

A direct consequence of the previous Lemma is the following 
\begin{Prop}
Let \(\gamma_t\) be a variation  of \(\gamma\) as described above.
\begin{enumerate}
 \item If \(\eta\in\Gamma(\gamma^\ast TM)\) then \eqref{first-variation-p-bienery} yields the equation for \(p\)-biharmonic curves
 \begin{align}
 \label{p-biharmonic-euler-lagrange}
   \nabla_{\gamma'}\nabla_{\gamma'}\big(|\nabla_{\gamma'}\gamma'|^{p-2}\nabla_{\gamma'}\gamma'\big)
 -|\nabla_{\gamma'}\gamma'|^{p-2}R^M(\gamma',\nabla_{\gamma'}\gamma')\gamma'=0.
 \end{align}
 \item If \(\eta\in\Gamma((\gamma^\ast TM)^\perp) \) then \eqref{first-variation-p-bienery} yields the equation for \(p\)-elastic curves
\begin{align}
\label{p-elastic-euler-lagrange}
  \bigg(( \nabla_{\gamma'}\nabla_{\gamma'}\big(|\nabla_{\gamma'}\gamma'|^{p-2}\nabla_{\gamma'}\gamma'\big)
 -|\nabla_{\gamma'}\gamma'|^{p-2}R^M(\gamma',\nabla_{\gamma'}\gamma')\gamma'\bigg)^\perp=0.
 \end{align}
\end{enumerate}
\end{Prop}

\subsection{The Frenet-equations}

For the further analysis it turns out to be useful to rewrite the equations for both
\(p\)-biharmonic \eqref{p-biharmonic-euler-lagrange} as well as \(p\)-elastic curves \eqref{p-elastic-euler-lagrange} in terms of its Frenet-frames.
To this end, we recall the following
\begin{Dfn}[Frenet-frame]
\label{dfn:frenet}
Let \(\gamma\colon I\to M\) be a curve which is parametrized with respect to arclength. Then its Frenet-frame is defined by the following set of
equations
\begin{align}
\label{frenet-frame}
F_1=&\gamma',\\
\nonumber\nabla_{\gamma'} F_1=&k_1F_2,\\
\nonumber\nabla_{\gamma'} F_i=&-k_{i-1}F_{i-1}+k_iF_{i+1},\qquad i=2,\ldots,n-1,\\
\nonumber\nonumber\vdots \\
\nonumber\nabla_{\gamma'} F_n=&-k_{n-1}F_{n-1}.
\end{align}
\end{Dfn}

For more details on Frenet-frames and their application in geometry
we refer to the book \cite{MR3676571}.

\begin{Lem}
Let \(\gamma\colon I\to M\) be a smooth curve 
parametrized by arclength
with its associated Frenet-frame.
Then the following equations hold
\begin{align*}
\nabla_{\gamma'}\nabla_{\gamma'}\big(|\nabla_{\gamma'}\gamma'|^{p-2}\nabla_{\gamma'}\gamma'\big)
=&\big((1-2p)k_1^{p-1}k_1'\big)F_1 \\
&+\big((p-1)(p-2)k_1^{p-3}k_1'^2+(p-1)k_1^{p-2}k_1''-k_1^{p+1}-k_1^{p-1}k_2^2\big)F_2 \\
&+\big(2(p-1)k_1^{p-2}k_1'k_2+k_1^{p-1}k_2'\big)F_3 \\
&+\big(k_1^{p-1}k_2k_3\big)F_4, \\
|\nabla_{\gamma'}\gamma'|^{p-2}R^M(\gamma',\nabla_{\gamma'}\gamma')\gamma'=&k_1^{p-1} R^M(F_1,F_2)F_1.
\end{align*}
\end{Lem}
\begin{proof}
This follows by a direct calculation using the Frenet-frame \eqref{frenet-frame}.
\end{proof}

\begin{Prop}
Let \(\gamma\colon I\to M\) be a smooth curve 
parametrized by arclength
with its associated Frenet-frame.
A curve is \(p\)-biharmonic, that is a critical point of \eqref{energy-p-biharmonic}, if the following system holds
\begin{align}
\label{p-biharmonic-curve}
0=&(1-2p)k_1^{p-1}k_1', \\
\nonumber 0=&(p-1)(p-2)k_1^{p-3}k_1'^2+(p-1)k_1^{p-2}k_1''-k_1^{p+1}-k_1^{p-1}k_2^2+k_1^{p-1}\langle R^M(F_1,F_2)F_2,F_1\rangle, \\
\nonumber 0=&2(p-1)k_1^{p-2}k_1'k_2+k_1^{p-1}k_2'-k_1^{p-1}\langle R^M(F_1,F_2)F_1,F_3\rangle,\\
\nonumber 0=&k_1^{p-1}k_2k_3-k_1^{p-1}\langle R^M(F_1,F_2)F_1,F_4\rangle,\\
\nonumber 0=&k_1^{p-1}\langle R^M(F_1,F_2)F_1,F_j\rangle,\qquad j=5,\ldots,n.
\end{align}
A curve is \(p\)-elastic if it solves the system \eqref{p-biharmonic-curve} neglecting the first equation.
\end{Prop}

\begin{proof}
This is a direct consequence of the previous Lemma.
\end{proof}

At this point we have to make the following case distinction:

\begin{Cor}
\begin{enumerate}
 \item If \(p\neq\frac{1}{2}\) the system \eqref{p-biharmonic-curve} reduces to the well-known system of equations that characterizes
 a proper biharmonic curve (see for example \cite{MR2250208}), i.e.
 \begin{equation}\label{biharmonic-curve}
\begin{cases}
k_1=const\neq 0, \\
\nonumber k_1^2+k_2^2=\langle R^M(F_1,F_2)F_2,F_1\rangle, \\
\nonumber k_2'=\langle R^M(F_1,F_2)F_1,F_3\rangle,\\
\nonumber k_2k_3=\langle R^M(F_1,F_2)F_1,F_4\rangle,\\
\nonumber \langle R^M(F_1,F_2)F_1,F_j\rangle =0,\qquad j=5,\ldots,n.
\end{cases}
\end{equation}
\item If \(p=\frac{1}{2}\) the first equation of \eqref{p-biharmonic-curve} is always satisfied and we obtain the system
\begin{equation}
\begin{cases}
\frac{3}{4}k_1^{-1}k_1'^2-\frac{1}{2}k_1''-k_1^3-k_1k_2^2+k_1\langle R^M(F_1,F_2)F_2,F_1\rangle=0, \\
\nonumber -k_1'k_2+k_1k_2'-k_1\langle R^M(F_1,F_2)F_1,F_3\rangle=0,\\
\nonumber k_2k_3-\langle R^M(F_1,F_2)F_1,F_4\rangle=0,\\
\nonumber \langle R^M(F_1,F_2)F_1,F_j\rangle =0,\qquad j=5,\ldots,n.
\end{cases}
\end{equation}
\end{enumerate}
\end{Cor}

Let us make the following remarks:

\begin{Bem}
\begin{enumerate}
\item Let us give some further explanations on the common features and differences
of the energy functionals \eqref{energy-p-biharmonic} and \eqref{p-elastic-energy}.
If we start with the energy for \(p\)-biharmonic curves \eqref{energy-p-biharmonic} and use the Frenet-equations \eqref{frenet-frame} then formally the energy functionals for \(p\)-biharmonic curves \eqref{energy-p-biharmonic} and 
for \(p\)-elastic curves \eqref{p-elastic-energy} coincide.
However, by employing the Frenet-equations \eqref{frenet-frame} in \eqref{energy-p-biharmonic} we tacitly assume that we are switching to an arclength parametrization
and thus, in general, both variational problems will be different.

\item In the case of \(p\)-biharmonic curves with \(p\neq\frac{1}{2}\) we always have to look for curves of constant
 geodesic curvature \(k_1\). This condition comes from the Frenet-equation that describes the tangential direction of the curve \(\gamma\).
 This is in sharp contrast to \(p\)-elastic curves which allow for solutions
 with non-constant curvature for all values of \(p\).
 However, for \(\frac{1}{2}\)-biharmonic curves the tangential equation from the Frenet-frame is automatically satisfied
 due to our special choice of \(p\) allowing for solutions of non-constant curvature.
 \item Note that also for \(p=\frac{1}{2}\) we always have the constant geodesic curvature solutions.
 However, there may be additional solutions and understanding the latter is the content of this article.
 \item We would like to point out that, in particular for \(p<2\), 
 one has to be very careful if a geodesic also is a solution of the
 equation for \(p\)-biharmonic curves as there might be a negative power
 of the geodesic curvature involved. One way to get around this problem is
 to remove geodesics from the spaces of curves we are working with.
\end{enumerate}
\end{Bem}

Motivated from the previous remark we make the following definition:
\begin{Dfn}
A curve \(\gamma\colon I\to M\) is called a proper \(\frac{1}{2}\)-biharmonic curve   
if it has non-constant geodesic curvature.
\end{Dfn}

\begin{Bem}
Note that the notion of proper \(\frac{1}{2}\)-biharmonic curve
also excludes the case of geodesics as these have constant geodesic curvature
zero.
\end{Bem}

\begin{Bem}
We can also show that \(p\)-biharmonic curves can have non-constant
geodesic curvature without using any Frenet-frame.
Testing \eqref{p-biharmonic-euler-lagrange} with \(\gamma'\) we get
\begin{align*}
0=\langle \nabla_{\gamma'}\nabla_{\gamma'}\big(|\nabla_{\gamma'}\gamma'|^{p-2}\nabla_{\gamma'}\gamma'\big),\gamma'\rangle,
\end{align*}
from which we immediately deduce
\begin{align*}
0=-2\frac{d}{ds}|\nabla_{\gamma'}\gamma'|^p
+|\nabla_{\gamma'}\gamma'|^{p-2}\langle\nabla_{\gamma'}\nabla_{\gamma'}\gamma',
\nabla_{\gamma'}\gamma'\rangle.
\end{align*}
Here, we used that \(\gamma\) is parametrized with respect to arclength. 
The above equation directly implies that 
\begin{align*}
0=\big(-2+\frac{1}{p}\big)\frac{d}{ds}|\nabla_{\gamma'}\gamma'|^p
\end{align*}
and holds trivially for \(p=\frac{1}{2}\).
\end{Bem}

Let us also make a comment on the Euler-Lagrange method for \(p\)-biharmonic curves
which is a straightforward generalization of the corresponding result for biharmonic curves from  
\cite{MR3045700}.

We define the Lagrangian for \(p\)-biharmonic curves by
\begin{align}
\label{lagrangian}
\mathcal{L}_{2,p}:=|\tau(\gamma)|^p=g\big(\tau(\gamma),\tau(\gamma)\big)^\frac{p}{2}.
\end{align}
Then, either by a direct calculation, or by the general theory
of one-dimensional variational problems, we obtain the following
\begin{Satz}
\label{theorem-euler-lagrange}
The equation for \(p\)-biharmonic curves
\begin{align*}
\nabla_{\gamma'}\nabla_{\gamma'}\big(|\nabla_{\gamma'}\gamma'|^{p-2}\nabla_{\gamma'}\gamma'\big)
-|\nabla_{\gamma'}\gamma'|^{p-2}R^M(\gamma',\nabla_{\gamma'}\gamma')\gamma'=0
\end{align*}
is equivalent to the system of Euler-Lagrange equations
\begin{align}
\label{euler-lagrange}
\frac{d^2}{ds^2}\big(\frac{\partial\mathcal{L}_{2,p}}{\partial\gamma''^k}\big)
-\frac{d}{ds}\big(\frac{\partial\mathcal{L}_{2,p}}{\partial\gamma'^k}\big)
+\frac{\partial\mathcal{L}_{2,p}}{\partial\gamma^k}=0,
\qquad k=1,\ldots,n,
\end{align}
where the Lagrangian \(\mathcal{L}_{2,p}\) is defined in \eqref{lagrangian}.
\end{Satz}

\subsection{The second variation formula}
In this subsection we calculate the second variation of the energy
of \(p\)-biharmonic curves \eqref{energy-p-biharmonic}.

Again, we fix a small number \(\epsilon>0\) and let \(\gamma_t\colon (-\epsilon,\epsilon)\times I\to M\) be a variation
of the curve \(\gamma\) satisfying \(\frac{\nabla\gamma_t}{\partial t}\big|_{t=0}=\eta\), where
\(\eta\in\Gamma(\gamma^\ast TM)\) and \(I\subset\R\) is a closed interval.

\begin{Prop}
The second variation of \eqref{energy-p-biharmonic} 
evaluated at a critical point, that is a solution of \eqref{p-biharmonic-euler-lagrange},
is given by the following expression
\begin{align}
\label{eqn:second-variation}
\frac{d^2}{dt^2}&\big|_{t=0}\frac{1}{p}\int_I|\nabla_{\gamma_t'}\gamma_t'|^pds \\
=&\nonumber
\int_I\langle\eta, R^M(\eta,\gamma')\nabla_{\gamma'}\bigg(|\nabla_{\gamma'}\gamma'|^{p-2}\nabla_{\gamma'}\gamma'\bigg)\rangle ds 
-2\int_I|\nabla_{\gamma'}\gamma'|^{p-2}\langle
R^M(\nabla_{\gamma'}\eta,\nabla_{\gamma'}\gamma')\gamma',\eta\rangle ds  \\
&\nonumber+(p-2)\int_I
|\nabla_{\gamma'}\gamma'|^{p-4}
|\langle\nabla_{\gamma'}\nabla_{\gamma'}\eta,\nabla_{\gamma'}\gamma'\rangle|^2
 ds \\
&\nonumber+2(p-2)\int_I
|\nabla_{\gamma'}\gamma'|^{p-4}
\langle\nabla_{\gamma'}\nabla_{\gamma'}\eta,\nabla_{\gamma'}\gamma'
\rangle 
\langle R^M(\eta,\gamma')\gamma',\nabla_{\gamma'}\gamma'\rangle ds\\
&\nonumber+\int_I|\nabla_{\gamma'}\gamma'|^{p-2}|\nabla_{\gamma'}\nabla_{\gamma'}\eta|^2ds
+2\int_I|\nabla_{\gamma'}\gamma'|^{p-2}\langle
R^M(\eta,\gamma')\gamma',\nabla_{\gamma'}\nabla_{\gamma'}\eta
\rangle ds
\\
&\nonumber+(p-2)\int_I|\nabla_{\gamma'}\gamma'|^{p-4} 
|\langle R^M(\gamma,\nabla_{\gamma'}\gamma')\gamma',\eta\rangle|^2 ds \\
&\nonumber-\int_I|\nabla_{\gamma'}\gamma'|^{p-2}
\langle(\nabla_\eta R^M)(\gamma',\nabla_{\gamma'}\gamma')\gamma',\eta\rangle ds 
+\int_I|\nabla_{\gamma'}\gamma'|^{p-2}
|R^M(\gamma',\eta)\gamma'|^2 ds \\
&\nonumber-\int_I|\nabla_{\gamma'}\gamma'|^{p-2}
\langle R^M(\gamma',\nabla_{\gamma'}\gamma')\nabla_{\gamma'}\eta,\eta\rangle ds.
\end{align}

\end{Prop}

\begin{proof}
Let us recall the first variation formula derived in Lemma 
\ref{first-variation-p-bienery}, that is
\begin{align*}
\frac{d}{dt}\big|_{t=0}\frac{1}{p}\int_I|\nabla_{\gamma_t'}\gamma_t'|^pds=
\int_I\big(\langle \eta,
\nabla_{\gamma'}\nabla_{\gamma'}\big(|\nabla_{\gamma'}\gamma'|^{p-2}\nabla_{\gamma'}\gamma'\big)
 -|\nabla_{\gamma'}\gamma'|^{p-2}R^M(\gamma',\nabla_{\gamma'}\gamma')\gamma'
\rangle\big) ds
\end{align*}
and the identity
\begin{align*}
\frac{\nabla}{\partial t}\big|_{t=0}\nabla_{\gamma_t'}\gamma_t'=\nabla_{\gamma'}\nabla_{\gamma'}\eta
+R^M(\eta,\gamma')\gamma'.
\end{align*}
Hence, we have
\begin{align*}
\frac{d}{dt}\big|_{t=0}|\nabla_{\gamma_t'}\gamma_t'|^{p-2}=
(p-2)|\nabla_{\gamma'}\gamma'|^{p-4}\big(
\langle\nabla_{\gamma'}\nabla_{\gamma'}\eta,\nabla_{\gamma'}\gamma'\rangle
+\langle R^M(\eta,\gamma')\gamma',\nabla_{\gamma'}\gamma'\rangle\big).
\end{align*}

By a direct calculation we find
\begin{align*}
\frac{\nabla}{\partial t}\big|_{t=0}&\bigg(\nabla_{\gamma_t'}\nabla_{\gamma_t'}\big(|\nabla_{\gamma_t'}\gamma_t'|^{p-2}\nabla_{\gamma_t'}\gamma_t'\big)\bigg)\\
=&R^M(\eta,\gamma')\nabla_{\gamma'}\bigg(|\nabla_{\gamma'}\gamma'|^{p-2}\nabla_{\gamma'}\gamma'\bigg)
+\nabla_{\gamma_t'}\bigg(\frac{\nabla}{\partial t}\big|_{t=0}\nabla_{\gamma_t'}\big(|\nabla_{\gamma_t'}\gamma_t'|^{p-2}\nabla_{\gamma_t'}\gamma_t'\big)\bigg) \\
=&R^M(\eta,\gamma')\nabla_{\gamma'}\bigg(|\nabla_{\gamma'}\gamma'|^{p-2}\nabla_{\gamma'}\gamma'\bigg)
+\nabla_{\gamma'}\bigg(|\nabla_{\gamma'}\gamma'|^{p-2}R^M(\eta,\gamma')\nabla_{\gamma'}\gamma'\bigg) \\
&+\nabla_{\gamma_t'}\nabla_{\gamma_t'}\bigg(\frac{\nabla}{\partial t}\big|_{t=0}\big(|\nabla_{\gamma_t'}\gamma_t'|^{p-2}\nabla_{\gamma_t'}\gamma_t'\big)\bigg) \\
=&R^M(\eta,\gamma')\nabla_{\gamma'}\bigg(|\nabla_{\gamma'}\gamma'|^{p-2}\nabla_{\gamma'}\gamma'\bigg)
+\nabla_{\gamma'}\bigg(|\nabla_{\gamma'}\gamma'|^{p-2}R^M(\eta,\gamma')\nabla_{\gamma'}\gamma'\bigg) \\
&+(p-2)\nabla_{\gamma'}\nabla_{\gamma'}\bigg(
|\nabla_{\gamma'}\gamma'|^{p-4}\big(
\langle\nabla_{\gamma'}\nabla_{\gamma'}\eta,\nabla_{\gamma'}\gamma'\rangle\nabla_{\gamma'}\gamma'
+\langle R^M(\eta,\gamma')\gamma',\nabla_{\gamma'}\gamma'\rangle\big)\nabla_{\gamma'}\gamma'\bigg)\\
&+\nabla_{\gamma'}\nabla_{\gamma'}\bigg(|\nabla_{\gamma'}\gamma'|^{p-2}
(\nabla_{\gamma'}\nabla_{\gamma'}\eta
+R^M(\eta,\gamma')\gamma')
\bigg).
\end{align*}

In addition, we find
\begin{align*}
\frac{\nabla}{\partial t}\big|_{t=0}&\bigg(|\nabla_{\gamma'_t}\gamma_t'|^{p-2}R^M(\gamma_t',\nabla_{\gamma_t'}\gamma_t')\gamma_t'\bigg) \\
=&(p-2)|\nabla_{\gamma'}\gamma'|^{p-4}\bigg(
\langle\nabla_{\gamma'}\nabla_{\gamma'}\eta,\nabla_{\gamma'}\gamma'\rangle
+\langle R^M(\eta,\gamma')\gamma',\nabla_{\gamma'}\gamma'\rangle\bigg)R^M(\gamma',\nabla_{\gamma'}\gamma')\gamma'\big) \\
&+|\nabla_{\gamma'}\gamma'|^{p-2}
\bigg((\nabla_\eta R^M)(\gamma',\nabla_{\gamma'}\gamma')\gamma'
+R^M(\nabla_{\gamma'}\eta,\nabla_{\gamma'}\gamma')\gamma' \\
&\hspace{2cm}+R^M(\gamma',R^M(\eta,\gamma')\gamma')\gamma' 
+R^M(\gamma',\nabla_{\gamma'}\nabla_{\gamma'}\eta)\gamma'
+R^M(\gamma',\nabla_{\gamma'}\gamma')\nabla_{\gamma'}\eta\bigg).
\end{align*}

Combining the previous equations we get
\begin{align}
\label{eqn:second-variation-a}
\frac{d^2}{dt^2}&\big|_{t=0}\frac{1}{p}\int_I|\nabla_{\gamma_t'}\gamma_t'|^pds \\
=&\nonumber
\int_I\langle\eta, R^M(\eta,\gamma')\nabla_{\gamma'}\bigg(|\nabla_{\gamma'}\gamma'|^{p-2}\nabla_{\gamma'}\gamma'\bigg)\rangle ds 
+\int_I\langle\eta,\nabla_{\gamma'}\bigg(|\nabla_{\gamma'}\gamma'|^{p-2}R^M(\eta,\gamma')\nabla_{\gamma'}\gamma'\bigg)\rangle ds \\
&\nonumber+(p-2)\int_I\langle\eta,\nabla_{\gamma'}\nabla_{\gamma'}\bigg(
|\nabla_{\gamma'}\gamma'|^{p-4}
\langle\nabla_{\gamma'}\nabla_{\gamma'}\eta,\nabla_{\gamma'}\gamma'\rangle\nabla_{\gamma'}\gamma'\bigg)\rangle ds \\
&\nonumber+(p-2)\int_I\langle\eta,\nabla_{\gamma'}\nabla_{\gamma'}\bigg(
|\nabla_{\gamma'}\gamma'|^{p-4}
\langle R^M(\eta,\gamma')\gamma',\nabla_{\gamma'}\gamma'\rangle\nabla_{\gamma'}\gamma'\bigg)
\rangle ds\\
&\nonumber+\int_I\langle\eta,\nabla_{\gamma'}\nabla_{\gamma'}\bigg(|\nabla_{\gamma'}\gamma'|^{p-2}
\nabla_{\gamma'}\nabla_{\gamma'}\eta\bigg)\rangle ds
+\int_I\langle\eta,\nabla_{\gamma'}\nabla_{\gamma'}\bigg(|\nabla_{\gamma'}\gamma'|^{p-2}
R^M(\eta,\gamma')\gamma'
\bigg)\rangle ds
\\
&\nonumber-(p-2)\int_I|\nabla_{\gamma'}\gamma'|^{p-4}
\langle\nabla_{\gamma'}\nabla_{\gamma'}\eta,\nabla_{\gamma'}\gamma'\rangle
\langle R^M(\gamma',\nabla_{\gamma'}\gamma')\gamma',\eta\rangle ds \\
&\nonumber-(p-2)\int_I|\nabla_{\gamma'}\gamma'|^{p-4} 
\langle R^M(\eta,\gamma')\gamma',\nabla_{\gamma'}\gamma'\rangle\langle R^M(\gamma',\nabla_{\gamma'}\gamma')\gamma',\eta\rangle ds \\
&\nonumber-\int_I|\nabla_{\gamma'}\gamma'|^{p-2}
\langle(\nabla_\eta R^M)(\gamma',\nabla_{\gamma'}\gamma')\gamma',\eta\rangle ds 
-\int_I|\nabla_{\gamma'}\gamma'|^{p-2}\langle
R^M(\nabla_{\gamma'}\eta,\nabla_{\gamma'}\gamma')\gamma',\eta\rangle ds \\
&\nonumber-\int_I|\nabla_{\gamma'}\gamma'|^{p-2}
\langle R^M(\gamma',R^M(\eta,\gamma')\gamma')\gamma',\eta\rangle ds 
-\int_I|\nabla_{\gamma'}\gamma'|^{p-2}
\langle R^M(\gamma',\nabla_{\gamma'}\nabla_{\gamma'}\eta)\gamma',\eta\rangle ds \\
&\nonumber
-\int_I|\nabla_{\gamma'}\gamma'|^{p-2}
\langle R^M(\gamma',\nabla_{\gamma'}\gamma')\nabla_{\gamma'}\eta,\eta\rangle ds\\
\nonumber:=&\sum_{i=1}^{13} J_i.
\end{align}

Now, we will manipulate the terms on the right-hand side
of \eqref{eqn:second-variation-a}.
First, we note that
\begin{align*}
J_2+J_{10}=-2\int_I|\nabla_{\gamma'}\gamma'|^{p-2}\langle
R^M(\nabla_{\gamma'}\eta,\nabla_{\gamma'}\gamma')\gamma',\eta\rangle ds.
\end{align*}

Concerning the \(J_3\)-term integration by parts yields
\begin{align*}
J_3=(p-2)\int_I|\nabla_{\gamma'}\gamma'|^{p-4}
|\langle\nabla_{\gamma'}\nabla_{\gamma'}\eta,
\nabla_{\gamma'}\gamma'\rangle|^2ds.
\end{align*}

Moreover, it is easy to see, using integration by parts again, that
\begin{align*}
J_4=(p-2)\int_I
\langle\nabla_{\gamma'}\nabla_{\gamma'}\eta,\nabla_{\gamma'}\gamma'
\rangle 
|\nabla_{\gamma'}\gamma'|^{p-4}
\langle R^M(\eta,\gamma')\gamma',\nabla_{\gamma'}\gamma'\rangle ds.
\end{align*}

Hence, we can conclude that 
\begin{align*}
J_4+J_7=2(p-2)\int_I
\langle\nabla_{\gamma'}\nabla_{\gamma'}\eta,\nabla_{\gamma'}\gamma'
\rangle 
|\nabla_{\gamma'}\gamma'|^{p-4}
\langle R^M(\eta,\gamma')\gamma',\nabla_{\gamma'}\gamma'\rangle ds.
\end{align*}

For the \(J_5\)-term we get
\begin{align*}
J_5
=\int_I|\nabla_{\gamma'}\gamma'|^{p-2}|\nabla_{\gamma'}\nabla_{\gamma'}\eta|^2ds.
\end{align*}

In addition, similar arguments as before show that
\begin{align*}
J_6+J_{12}=
2\int_I|\nabla_{\gamma'}\gamma'|^{p-2}
\langle R^M(\eta,\gamma')\gamma',\nabla_{\gamma'}\nabla_{\gamma'}\eta\rangle ds.
\end{align*}

Regarding the \(J_{8}\)-term, 
we use the symmetries of the Riemannian curvature tensor and find
\begin{align*}
\langle R^M(\eta,\gamma')\gamma',\nabla_{\gamma'}\gamma'\rangle\langle R^M(\gamma',\nabla_{\gamma'}\gamma')\gamma',\eta\rangle
=-|\langle R^M(\gamma,\nabla_{\gamma'}\gamma')\gamma',\eta\rangle|^2.
\end{align*}

Concerning the \(J_{11}\)-term, we use
\begin{align*}
\langle R^M(\gamma',R^M(\eta,\gamma')\gamma')\gamma',\eta\rangle
=-|R^M(\gamma',\eta)\gamma'|^2.
\end{align*}
Inserting these identities into \eqref{eqn:second-variation-a}
completes the proof.
\end{proof}

\begin{Bem}
Note that for \(p=2\) equation \eqref{eqn:second-variation}
coincides with the calculation of Jiang \cite{MR886529} for the second
variation of the bienergy for maps between arbitrary Riemannian manifolds.
\end{Bem}

In the following we will assume that the manifold \(M\) is a space form,
that is it has a metric of constant curvature \(K\).
In this case the Riemannian curvature tensor can be written in the
following form
\begin{align}
\label{curvature-tensor-constant}
R^{M}(X,Y)Z=K\big(\langle Z,Y\rangle X-\langle Z,X\rangle Y\big)
\end{align}
for vector fields \(X,Y,Z\).

\begin{Prop}
\label{prop-second-variation-spaceform}
Let \(\gamma\colon I\to M\) be a \(p\)-biharmonic curve and 
suppose that \(M\) is a space form of constant curvature \(K\).
Then, the second variation of \eqref{energy-p-biharmonic} 
evaluated at a critical point, that is a solution of 
\eqref{p-biharmonic-euler-lagrange},
simplifies as follows
\begin{align}
\label{eqn:second-variation-space-form}
\frac{d^2}{dt^2}&\big|_{t=0}\frac{1}{p}\int_I|\nabla_{\gamma_t'}\gamma_t'|^pds \\
=&\nonumber(p-2)\int_I
|\nabla_{\gamma'}\gamma'|^{p-4}
|\langle\nabla_{\gamma'}\nabla_{\gamma'}\eta,\nabla_{\gamma'}\gamma'\rangle|^2
 ds +\int_I|\nabla_{\gamma'}\gamma'|^{p-2}|\nabla_{\gamma'}\nabla_{\gamma'}\eta|^2ds \\
&\nonumber
+(p-2)K^2\int_I|\nabla_{\gamma'}\gamma'|^{p-4} 
\big|\langle\eta,\gamma'\rangle\langle\gamma',\nabla_{\gamma'}\gamma'\rangle-|\gamma'|^2\langle\nabla_{\gamma'}\gamma',\eta\rangle\big|^2 ds \\ 
&\nonumber+K^2\int_I|\nabla_{\gamma'}\gamma'|^{p-2}
\big|\langle\eta,\gamma'\rangle\gamma'-|\gamma'|^2\eta\big|^2 ds\\
&\nonumber+K\int_I
\bigg(|\eta|^2\langle\gamma',\nabla_{\gamma'}\big(|\nabla_{\gamma'}\gamma'|^{p-2}\nabla_{\gamma'}\gamma'\big)\rangle
-\langle\eta,\gamma'\rangle\langle\eta,\nabla_{\gamma'}\big(|\nabla_{\gamma'}\gamma'|^{p-2}\nabla_{\gamma'}\gamma'\big)\rangle
\bigg) ds \\
&\nonumber
-2K\int_I|\nabla_{\gamma'}\gamma'|^{p-2}
\big(\langle\gamma',\nabla_{\gamma'}\gamma'\rangle
\langle\nabla_{\gamma'}\eta,\eta\rangle
-\langle\nabla_{\gamma'}\eta,\gamma'\rangle
\langle\nabla_{\gamma'}\gamma',\eta\rangle\big) ds  \\
&\nonumber
+2(p-2)K\int_I
\langle\nabla_{\gamma'}\nabla_{\gamma'}\eta,\nabla_{\gamma'}\gamma'
\rangle 
|\nabla_{\gamma'}\gamma'|^{p-4}
\big(|\gamma'|^2
\langle\eta,\nabla_{\gamma'}\gamma'\rangle
-\langle\eta,\gamma'\rangle\langle\gamma',\nabla_{\gamma'}\gamma'\rangle
\big) ds\\
&\nonumber-2K\int_I|\nabla_{\gamma'}\gamma'|^{p-2}
\big(\langle\gamma',\nabla_{\gamma'}\nabla_{\gamma'}\eta\rangle
\langle\gamma',\eta\rangle
-|\gamma'|^2\langle\nabla_{\gamma'}\nabla_{\gamma'}\eta,\eta\rangle\big)
ds \\
&\nonumber
-K\int_I|\nabla_{\gamma'}\gamma'|^{p-2}
\big(\langle\nabla_{\gamma'}\gamma',\nabla_{\gamma'}\eta\rangle
\langle\gamma',\eta\rangle
-\langle\nabla_{\gamma'}\eta,\gamma'\rangle\langle\nabla_{\gamma'}\gamma',\eta\rangle\big) ds.
\end{align}
\end{Prop}
\begin{proof}
The claim follows by inserting the definition of the curvature tensor of a space form \eqref{curvature-tensor-constant} into \eqref{eqn:second-variation}
and a direct calculation. Note that the term proportional to 
\(\nabla R\) drops out due to the assumption of constant curvature.
\end{proof}

\subsection{Qualitative aspects of \texorpdfstring{$p$}{}-biharmonic curves}
It is well-known that every biharmonic curve on a manifold with negative curvature
must be a geodesic. In the case of \(p\)-biharmonic maps we have a similar result,
see for example \cite{MR3393120}.
Here, we prove the following 

\begin{Satz}
\label{theorem-liouville}
Let \(\gamma\colon I\to M\) be a \(p\)-biharmonic curve parametrized by arclength.
If \(M\) has non-positive curvature then \(\gamma\) must be a geodesic.
\end{Satz}

\begin{proof}
From the second equation of the system \eqref{p-biharmonic-curve} we get
\begin{align*}
k_1^2+k_2^2=\langle R^M(F_1,F_2)F_2,F_1\rangle
\end{align*}
which immediately completes the proof.
\end{proof}

\begin{Bem}
We can also show the above claim by a different method in the case
that \(I\subset\R\) is a closed interval.
We test the equation for \(p\)-biharmonic curves \eqref{p-biharmonic-euler-lagrange} with 
\(|\nabla_{\gamma'}\gamma'|^{p-2}\nabla_{\gamma'}\gamma'\) 
(assuming that \(\gamma\) is non-geodesic)
and obtain 
\begin{align*}
\int_I\langle\nabla_{\gamma'}\nabla_{\gamma'}
\big(|\nabla_{\gamma'}\gamma'|^{p-2}\nabla_{\gamma'}\gamma'\big),
|\nabla_{\gamma'}\gamma'|^{p-2}\nabla_{\gamma'}\gamma'\rangle ds
+\int_I|\nabla_{\gamma'}\gamma'|^{2p-4}
\langle R^M(\gamma',\nabla_{\gamma'}\gamma')\nabla_{\gamma'}\gamma',\gamma'\rangle ds
=0.
\end{align*}
Using the assumption of \(M\) having negative sectional curvature 
together with integration by parts we can deduce that
\begin{align*}
\int_I|\nabla_{\gamma'}
\big(|\nabla_{\gamma'}\gamma'|^{p-2}\nabla_{\gamma'}\gamma'\big)|^2ds
\leq 0
\end{align*}
which implies that \(\nabla_{\gamma'}
\big(|\nabla_{\gamma'}\gamma'|^{p-2}\nabla_{\gamma'}\gamma'\big)=0\). 
Testing this equation with \(\gamma'\) and using integration by parts
once more we find
\begin{align*}
0=\int_I|\nabla_{\gamma'}\gamma'|^pds
\end{align*}
from which we deduce \(\nabla_{\gamma'}\gamma'=0\) completing the proof.
\end{Bem}

\begin{Bem}
We would like to point out once more that one has to be very careful
with results such as Theorem
\ref{theorem-liouville} since there might be negative 
powers of the geodesic curvature involved.
\end{Bem}

\section{\texorpdfstring{$p$}{}-biharmonic curves in two dimensions}
In this section we study \(p\)-biharmonic curves on two-dimensional domains.
We will pay special attention to the Euclidean plane \(\R^2\) modelling the 
zero curvature case \(K=0\), the round sphere \(\s^2\) with constant curvature \(K=1\)
and hyperbolic surfaces with constant curvature \(K=-1\).

We have already seen that for \(p\neq\frac{1}{2}\) the solutions of the equation for 
\(p\)-biharmonic curves are precisely biharmonic curves, see \cite{MR1884940}, 
such that we do not have
to further investigate this case.

On the other hand, the case of \(p=\frac{1}{2}\)
has already been studied in the context of \(p\)-elastic curves on surfaces 
and critical points of \eqref{energy-blaschke} in
\cite{MR2031956,MR2007599,MR4028018,MR1548296}. Our results are very close
to the ones obtained in that references.

On a two-dimensional manifold we have the Frenet-frame \(\{T,N\}\) satisfying the following equations
\begin{align}
\label{frenet-surface}
\nabla_TT=kN,\qquad \nabla_TN=-kT
\end{align}
corresponding to \eqref{frenet-frame} where we relabelled \(k_1=k,F_1=T,F_2=N\)
 in order to be consistent with the notation usually 
employed in the case of a curve on a surface.

Writing the equation for \(p\)-biharmonic curves in terms of the Frenet-frame \eqref{frenet-surface}
we obtain the system
\begin{align*}
(p-1)(p-2)k ^{-1}k '^2N+(p-1)k ''N-(2p-1)k k 'T-k ^3N-k  R^M(T,N)T=0.
\end{align*}
Taking the scalar product with \(T\) and \(N\) we get the two equations
\begin{align*}
(2p-1)k k '&=0,\\
\nonumber(p-1)(p-2)k ^{-1}k '^2+(p-1)k ''-k ^3+k  K&=0.
\end{align*}

In the case that \(p\neq\frac{1}{2}\) we get the well-known system
\begin{align}
k =const,\qquad k^2=K
\end{align}
which characterizes biharmonic curves on surfaces \cite{MR1884940}.

However, in the case that \(p=\frac{1}{2}\) the first equation is satisfied trivially
and we get the following

\begin{Prop}
Let \(\gamma\colon I\to M\) be a curve parametrized by arclength and \(M\) a surface with Gauss curvature \(K\).
Then the curve \(\gamma\) is proper \(\frac{1}{2}\)-biharmonic if its geodesic curvature \(k\)
is a non-constant solution of the following ordinary differential equation
\begin{align}
\label{ode-k-surface}
\frac{3}{4}k'^2-\frac{1}{2}k''k-k^4+k^2K=0.
\end{align}
\end{Prop}
It can be directly seen that \(k^2=K=const\) gives a solution of \eqref{ode-k-surface}.
Hence, in the following we will assume that \(k\neq const\), 
then the following identity holds
\begin{align}
\label{identity-k}
k ^\frac{5}{2}\big(k ^{-\frac{1}{2}}\big)''=\frac{3}{4}k '^2-\frac{1}{2}k ''k 
\end{align}
and \eqref{ode-k-surface} simplifies to
\begin{align*}
\big(k ^{-\frac{1}{2}}\big)''-k^\frac{3}{2}+k ^{-\frac{1}{2}}K=0.
\end{align*}
Setting \(f:=k ^{-\frac{1}{2}}\), which we assume to be non-constant,
and multiplying by \(f'\) we find 
\begin{align*}
0=(f'^2)'+K(f^2)'+\big(\frac{1}{f^2}\big)'.
\end{align*}
This allows us to obtain the simpler equation
\begin{align*}
c_1=f'^2+Kf^2+\frac{1}{f^2}
\end{align*}
for some integration constant \(c_1\in\R\).

To further simplify this equation we define \(h(s):=\frac{1}{2}f^2(s)\)
such that we are left with the differential equation
\begin{align}
\label{ode-g}
0=h'^2+4Kh^2-2c_1h+1.
\end{align}
The solutions of \eqref{ode-g} then give the geodesic curvature via the assignment \(k(s)=\frac{1}{2h(s)}\).
The next theorem discusses the solutions of \eqref{ode-g} in the three cases \(K=1,0,-1\) which is very similar to the mathematical analysis carried out in  \cite{MR2031956,MR2007599,MR4028018,MR1548296} in the context
of generalized elastic curves and is included in \cite[Proposition 3.1]{MR3774309}.

\begin{Satz}
\label{thm:curvature-surface}
Let \(\gamma\colon I\to\ M\) be a smooth curve parametrized by arclength and 
\(M\) a surface with constant curvature \(K\).
\begin{enumerate}
 \item If \(M=\R^2\) with the flat metric the curve \(\gamma\) is proper \(\frac{1}{2}\)-biharmonic if its geodesic curvature
is given by
\begin{align}
k(s)=\frac{c_1}{c_1^2(c_2+s)^2+1},
\end{align}
where \(c_1,c_2\in\R\).
 \item If  \(M=\s^2\) with the round metric with curvature \(K=1\)
 the curve \(\gamma\) is proper \(\frac{1}{2}\)-biharmonic if its geodesic curvature is given by
\begin{align}
k(s)=\frac{2}{c_1+\sqrt{c_1^2-4}\sin 2(c_3+s)},
\end{align}
where \(c_1,c_3\in\R\).
\item If \(M\) is a hyperbolic surface with a metric of constant curvature \(K=-1\)
the curve \(\gamma\) is proper \(\frac{1}{2}\)-biharmonic if its geodesic curvature is given by
\begin{align}
k(s)=\frac{4e^{2(c_4+s)}}{(c_1-e^{2(c_4+s)})^2+4},
\end{align}
where \(c_1,c_4\in\R\).
\end{enumerate}
\end{Satz}

\begin{Bem}
Note that by shifting the curve parameter \(s\) appropriately
the constants \(c_2,c_3,c_4\) can assumed to be zero.
\end{Bem}

\begin{proof}
In the case of \(M=\R^2\) with the flat Euclidean metric we have \(K=0\) and \eqref{ode-g} simplifies to
\begin{align*}
0=h'^2-2c_1h+1.
\end{align*}
This equation can be integrated as
\begin{align*}
h(s)=\frac{1+c_1^2(c_2+s)^2}{2c_1},
\end{align*}
where the integration constant \(c_2\in\R\) is determined via the initial data.

In the case of \(M=\s^2\) with the round metric with \(K=1\) equation \eqref{ode-g}
acquires the form
\begin{align*}
0=h'^2+4h^2-2c_1h+1.
\end{align*}
This equation can again be integrated directly as
\begin{align*}
h(s)=\frac{c_1}{4}+\frac{1}{4}\sqrt{c_1^2-1}\sin 2(c_3+s),
\end{align*}
where, again, the integration constant \(c_3\in\R\) is determined by the initial data.

On a hyperbolic surface we can always choose a metric of constant Gaussian curvature \(K=-1\) such that \eqref{ode-g} reduces to
\begin{align*}
0=h'^2-4h^2-2c_1h+1.
\end{align*}

A solution of this equation is given by
\begin{align*}
h(s)=\frac{1}{8}e^{-2(c_4+s)}\big((c_1-e^{2(c_4+s)})^2+4\big),
\end{align*}
with the integration constant \(c_4\in\R\).
Now, the proof is complete.
\end{proof}

\subsection{Existence results via magnetic geodesics}
In the case of \(M=\s^2,\) or \(M\) being a hyperbolic surface we can make use
of the results obtained for (closed) \emph{magnetic geodesics} by Schneider
\cite{MR2788659,MR2959932} to obtain an existence result for proper
\(\frac{1}{2}\)-biharmonic curves.

Magnetic geodesics describe the trajectory of a particle evolving
in a Riemannian manifold \(M\) under the influence of an external magnetic field.
However, there also is a more geometric point of view on such curves
by asking: Given a function \(k\colon M\to\R\), does there exist a (closed)
curve \(\gamma\) on \(M\) with geodesic curvature \(k\)?
In order to answer this question, or to find a magnetic geodesic,
one has to solve the equation
\begin{align}
\label{eq:magnetic}
\nabla_{\gamma'}\gamma'=k J_g(\gamma'),
\end{align}
where \(J_g\) represents the rotation by \(\pi/2\) in \(TM\) measured
with respect to the Riemannian metric \(g\), which corresponds
to the Frenet-equations \eqref{frenet-surface}.
Note that in the mathematics literature on magnetic geodesics one usually thinks of magnetic geodesics as closed curves
without explicitly mentioning this.

Now, we recall the following result (see \cite[Theorem 1.3]{MR2788659})
\begin{Satz}[Schneider, 2011]
\label{thm:schneider-s2}
Consider \(\s^2\) with the round metric of constant curvature \(K=1\).
For any positive constant \(k_0>0\) there exists a positive smooth function \(k\colon\s^2\to\R\), which can be chosen arbitrarily close to \(k_0\)
such that there are exactly two simple solutions of \eqref{eq:magnetic}.
\end{Satz}

In other words, this means that on the sphere with the round metric we can solve equation \eqref{eq:magnetic} for any prescribed curvature \(k\).

Recall that due to Theorem \ref{thm:curvature-surface}, case (2), the
geodesic curvature of a proper \(\frac{1}{2}\)-biharmonic curve on \(\s^2\)
with the round metric is given by
\begin{align*}
k(s)=\frac{2}{c_1+\sqrt{c_1^2-4}\sin 2s},
\end{align*}
where \(c_1\in\R\). 

Hence, if the integration constant \(c_1\) satisfies \(c_1>2\)
we have \(k(s)>0\) such that
we can apply Theorem \ref{thm:schneider-s2} to obtain
an existence result for proper \(\frac{1}{2}\)-biharmonic curves on \(\s^2\).

In addition, we recall the following result (see \cite[Theorem 1.2]{MR2959932}):
\begin{Satz}[Schneider, 2012]
\label{thm:schneider-hyperbolic}
Let \((M,g)\) be a closed oriented surface with negative Euler characteristic
\(\chi(M)\) and \(k\colon M\to\R\) be a positive function.
Assume that there exists \(K_0>0\) such that \(k\) and the Gaussian curvature \(K_g\)
of \((M,g)\) satisfy
\begin{align}
\label{eq:schneider-hyperbolic-ass}
k>(K_0)^\frac{1}{2}\quad\textrm{and} \quad K_g\geq -K_0.
\end{align}
Then, there exists an Alexandrov embedded curve \(\gamma\in C^2(S^1,M)\)
that solves \eqref{eq:magnetic} and the number of solutions is at least
\(-\chi(M)\) provided they are all non-degenerate.
\end{Satz}

Recall that due to Theorem \ref{thm:curvature-surface}, case (3), the
geodesic curvature of a proper \(\frac{1}{2}\)-biharmonic curve
on a hyperbolic surface with constant curvature \(K=-1\)
is given by
\begin{align*}
k(s)=\frac{4e^{2s}}{(c_1-e^{2s})^2+4},
\end{align*}
where \(c_1\in\R\). 
In order to satisfy the conditions \eqref{eq:schneider-hyperbolic-ass}
we could assume that \(s\geq 0\) or otherwise
we would need to require that \(c_1^2\) is sufficiently large.
Then, we can directly apply Theorem \ref{thm:schneider-hyperbolic} to get
an existence result for proper \(\frac{1}{2}\)-biharmonic curves
on hyperbolic surfaces.

\begin{Bem}
In \cite{MR2788659,MR2959932} there are also existence results
for magnetic geodesics on \(\s^2\) and hyperbolic surfaces
in case that these do not carry the standard metric of constant curvature.

In the case of constant Gauss curvature \(K\) there exist many results 
on \(\frac{1}{2}\)-elastic curves we could refer to \cite{MR2031956, MR2007599,MR2078687,MR4028018,MR4663707} in order to obtain the curvature
for proper \(\frac{1}{2}\)-biharmonic curves.
Note that these results rely on explicitly solving the equations for the curvature of the curve \eqref{ode-k-surface}.
However, in the case of a surface with non-constant Gauss curvature \(K\) it might be very difficult
to explicitly solve \eqref{ode-k-surface}
in which case the corresponding existence results from \cite{MR2788659,MR2959932}
for magnetic geodesics can provide solutions.
\end{Bem}

\begin{Bem}
As we have an explicit expression for the geodesic curvature of
\(\frac{1}{2}\)-biharmonic curves on surfaces by Theorem \ref{thm:curvature-surface}
it should also be possible to apply the heat flow for magnetic geodesics
studied in \cite{MR3778117} to get some existence results
for proper \(\frac{1}{2}\)-biharmonic curves.
\end{Bem}

\subsection{The stability of \texorpdfstring{$p$}{}-biharmonic curves in two dimensions}
In this subsection we investigate the stability of non-geodesic \(p\)-biharmonic curves
in two-dimensions.

In the following we will use the following notation
\begin{align*}
\operatorname{Hess}(E_{2,p}(\gamma))(\eta,\eta):=
\frac{d^2}{dt^2}&\big|_{t=0}E_{2,p}(\gamma_t),
\end{align*}
where \(\eta:=\frac{\nabla\gamma_t}{\partial t}\big|_{t=0}\in\Gamma(\gamma^\ast TM)\) denotes the variational vector field.

\begin{Prop}
Let \(\gamma\colon I\to M\) be a \(p\)-biharmonic curve parametrized by arclength
and suppose that \(M\) is a closed surface of constant curvature \(K\).
Then the Hessian of \eqref{energy-p-biharmonic} evaluated
at the unit normal \(N\) of the surface is given by
\begin{align}
\label{eq:second-variation-general}
\operatorname{Hess}(E_{2,p}(\gamma))(N,N)=
\int_I\big(
(p-1)k^{p+2}+k'^2k^{p-2}+(p-1)K^2k^{p-2}-(2p+2)Kk^p
\big)ds,
\end{align}
where \(k\) is the geodesic curvature of the curve \(\gamma\).
\end{Prop}

\begin{proof}
By the uniformization theorem every closed surface admits a metric 
of constant curvature \(K\). Hence, we can apply the second variation 
formula derived in Proposition \ref{prop-second-variation-spaceform}
and choose the unit normal \(N\) as variational vector field.

Note that by choosing \(\eta=N\) we get \(\langle\gamma',N\rangle=0\)
and also \(\langle\nabla_{\gamma'}N,N\rangle=0\).
Moreover, we also assume that \(\gamma\) is parametrized by arclength
such that we can use the Frenet-equations \eqref{frenet-surface}.

Inserting this choice into \eqref{eqn:second-variation-space-form}
we find
\begin{align*}
\operatorname{Hess}(E_{2,p}(\gamma))(N,N) 
=&\nonumber(p-2)\int_I
|\nabla_{\gamma'}\gamma'|^{p-4}
|\langle\nabla_{\gamma'}\nabla_{\gamma'} N,\nabla_{\gamma'}\gamma'\rangle|^2
 ds \\
&\nonumber
 +\int_I|\nabla_{\gamma'}\gamma'|^{p-2}|\nabla_{\gamma'}\nabla_{\gamma'}N|^2ds 
+(p-2)K^2\int_I|\nabla_{\gamma'}\gamma'|^{p-4} 
|\langle\nabla_{\gamma'}\gamma',N\rangle|^2 ds \\ 
&+K^2\int_I|\nabla_{\gamma'}\gamma'|^{p-2}ds
+K\int_I
\langle\gamma',\nabla_{\gamma'}\big(|\nabla_{\gamma'}\gamma'|^{p-2}\nabla_{\gamma'}\gamma'\big)\rangle
ds \\
&\nonumber
+3K\int_I|\nabla_{\gamma'}\gamma'|^{p-2}
\langle\nabla_{\gamma'}N,\gamma'\rangle
\langle\nabla_{\gamma'}\gamma',N\rangle ds  \\
&\nonumber
+2(p-2)K\int_I
\langle\nabla_{\gamma'}\nabla_{\gamma'}N,\nabla_{\gamma'}\gamma'
\rangle 
|\nabla_{\gamma'}\gamma'|^{p-4}
\langle N,\nabla_{\gamma'}\gamma'\rangle
 ds\\
&\nonumber+2K\int_I|\nabla_{\gamma'}\gamma'|^{p-2}
\langle\nabla_{\gamma'}\nabla_{\gamma'}N,N\rangle
ds.
\end{align*}

Let us manipulate the terms of the above equation which do not have a sign 
with the help of the Frenet-frame \eqref{frenet-surface} as
\begin{align*}
\int_I
\langle\gamma',\nabla_{\gamma'}\big(|\nabla_{\gamma'}\gamma'|^{p-2}\nabla_{\gamma'}\gamma'\big)\rangle
ds =&-\int_Ik^pds,  \\
\int_I|\nabla_{\gamma'}\gamma'|^{p-2}
\langle\nabla_{\gamma'}N,\gamma'\rangle
\langle\nabla_{\gamma'}\gamma',N\rangle ds  
=&-\int_Ik^pds, \\
\int_I
\langle\nabla_{\gamma'}\nabla_{\gamma'}N,\nabla_{\gamma'}\gamma'
\rangle 
|\nabla_{\gamma'}\gamma'|^{p-4}
\langle N,\nabla_{\gamma'}\gamma'\rangle
 ds=&-\int_Ik^pds,\\
\int_I|\nabla_{\gamma'}\gamma'|^{p-2}
\langle\nabla_{\gamma'}\nabla_{\gamma'}N,N\rangle
ds=&-\int_Ik^pds.
\end{align*}

Hence, we get
\begin{align*}
\operatorname{Hess}(E_{2,p}(\gamma))(N,N)
=&(p-2)\int_I
|\nabla_{\gamma'}\gamma'|^{p-4}
|\langle\nabla_{\gamma'}\nabla_{\gamma'}N,\nabla_{\gamma'}\gamma'\rangle|^2
 ds \\
 &+\int_I|\nabla_{\gamma'}\gamma'|^{p-2}|\nabla_{\gamma'}\nabla_{\gamma'}N|^2ds \\
&+(p-1)K^2\int_Ik^{p-2}ds 
-(2p+2)K\int_Ik^pds\\
=&(p-2)\int_Ik^{p+2}ds
+\int_I(k^{p+2}+k'^2k^{p-2})ds \\
&+(p-1)K^2\int_Ik^{p-2}ds 
-(2p+2)K\int_Ik^pds\\
=&\int_I\big(
(p-1)k^{p+2}+k'^2k^{p-2}+(p-1)K^2k^{p-2}-(2p+2)Kk^p
\big)ds,
\end{align*}
where we again used the Frenet-equations \eqref{frenet-surface}
in order to complete the proof.
\end{proof}

\begin{Bem}
We want to point out that the formula for the second variation \eqref{eq:second-variation-general} is different from the formula for the second variation of the elastic energy, which is \eqref{p-elastic-energy} for \(p=2\), obtained by   Langer and Singer \cite{MR772124} and also different from the stability of \(\frac{1}{2}\)-elastic curves, see \cite{MR4663707} and \cite{MR2031956}. While this may seem strange at first glance, there is a straightforward explanation for this phenomena: In the case of 
\(p\)-biharmonic curves one calculates the second variation of the \(L^p\)-norm of 
\(\nabla_{\gamma'}\gamma'\) and inserts the normal \(N\) in the quadratic form associated with the second variation, while in the case of \(p\)-elastic curves
one starts with the \(L^p\)-norm of the curvature of the curve, calculates the second variation and evaluates the associated quadratic form on \(N\). As these procedures are both slightly different in nature, one gets different results in the end.
\end{Bem}

\begin{Satz}
Let \(\gamma\colon I\to M\) be a \(p\)-biharmonic curve 
parametrized by arclength
and suppose that \(M\) is a closed surface of constant curvature \(K\).
Consider the second variation of \eqref{energy-p-biharmonic}
evaluated at a critical point
for a variation in the normal direction. Then the following
statements hold:
\begin{enumerate}
\item For \(p\neq\frac{1}{2}\)
  we get 
  \begin{align}
  \label{hessian-p-general}
\operatorname{Hess}(E_{2,p}(\gamma))(N,N) 
=-4K|I|k^{p}.
 \end{align}
\item For \(p=\frac{1}{2}\) we get
\begin{align}
\label{hessian-p12}
\operatorname{Hess}(E_{2,p}(\gamma))(N,N)
=\int_I\big(
-\frac{1}{2}k^\frac{5}{2}+k'^2k^{-\frac{3}{2}}
-\frac{1}{2}K^2k^{-\frac{3}{2}}-3Kk^\frac{1}{2}
\big)ds.
\end{align}

 \end{enumerate}
\end{Satz}
\begin{proof}
First, we consider the case \(p\neq\frac{1}{2}\).
Under this assumption all \(p\)-biharmonic curves on surfaces
are characterized by \(k^2=K=const\). Inserting 
into \eqref{eq:second-variation-general} we obtain the claim
after a direct calculation.

For \(p=\frac{1}{2}\) we get the statement directly from \eqref{eq:second-variation-general}.
\end{proof}

Let us make the following remarks concerning
the previous Theorem:

\begin{Bem}
\begin{enumerate}
 \item From \eqref{hessian-p-general} we can follow that non-geodesic 
 \(p\)-biharmonic curves for \(p\neq\frac{1}{2}\) on a surface with
 constant, positive curvature will always be unstable.
 \item In the case of \(p=\frac{1}{2}\) we can insert the solutions
  obtained in Theorem \ref{thm:curvature-surface} into equation
 \eqref{hessian-p12} in order to check the sign of the second variation.
 Unfortunately, this leads to a number of complicated expressions
 and it does not seem possible to give a definite answer.
 \item The instability of \(p\)-elastic curves on \(\s^2\) was recently
 shown in \cite{MR4663707} and one should thus expect that proper
 \(\frac{1}{2}\)-biharmonic curves on \(\s^2\) will be unstable as well.
\end{enumerate}

\end{Bem}

\begin{Bem}
\label{rem-yelin}
In \cite[Theorem 2.5]{MR4110268} it was shown that 
for a biharmonic map \(\phi\) from a closed Riemannian manifold \(M\) to a space form of positive sectional curvature \(K\) with \(|\tau(\phi)|^2=const\)
the second variation of the bienergy satisfies
\begin{align*}
\operatorname{Hess}(E_2(\phi))(\tau(\phi),\tau(\phi))=-4K|\tau(\phi)|^4\vol(M)\leq 0
\end{align*}
which clearly shows that such biharmonic maps are unstable.

This result is in perfect agreement with \eqref{hessian-p-general}:
If we choose \(p=2\), the variational vector field \(\eta=\tau(\gamma)=kN\)
and the equation for biharmonic curves on surfaces \(k^2=K=const\)
then we would get the same result.
\end{Bem}

\section{\texorpdfstring{$p$}{}-biharmonic curves on three-dimensional manifolds}
In this section we extend the previous analysis to the case of a three-dimensional
manifold, again assuming that we have a metric of constant curvature.

We choose the standard Frenet-frame \(\{T,N,B\}\) along the curve \(\gamma\)
which satisfies the following Frenet-equations
\begin{align}
\label{eq:frenet-3d}
\nabla_TT=kN,\qquad \nabla_TN=-kT+\tau B,\qquad \nabla_TB=-\tau N.
\end{align}
Here, \(k\) represents the \emph{geodesic curvature} of the curve \(\gamma\)
and \(\tau\) denotes the \emph{torsion} of the curve. Note that we again
relabelled the Frenet-frame \eqref{frenet-frame} in order to be consistent
with the notation usually applied in the three-dimensional setting.

Using the above Frenet-equations the condition for a curve to be \(p\)-biharmonic
is given by 
\begin{align}
0=&\big((p-1)(p-2)k ^{p-3}k '^2+(p-1)k ^{p-2}k ''-k ^{p+1}-k ^{p-1}\tau ^2\big)N \\
\nonumber&+\big((1-2p)k ^{p-1}k '\big)T \\
\nonumber&+\big(2(p-1)k ^{p-2}k '\tau +k ^{p-1}\tau '\big)B \\
\nonumber&-k ^{p-1}R^M(T,N)T.
\end{align}

This directly implies the following
\begin{Prop}
Let \(\gamma\colon I\to M\), where \(M\) is a three-dimensional Riemannian manifold, be a curve which is parametrized with respect to arclength together with
its associated Frenet-frame \(\{T,N,B\}\). Then \(\gamma\) is \(p\)-biharmonic if
the following equations hold
\begin{align*}
0=&(1-2p)k ^{p-1}k ',\\
0=&(p-1)(p-2)k ^{p-3}k '^2+(p-1)k ^{p-2}k ''-k ^{p+1}-k ^{p-1}\tau ^2+k ^{p-1}\langle R^M(T,N)N,T\rangle,\\
0=&2(p-1)k ^{p-2}k '\tau +k ^{p-1}\tau '-k ^{p-1}\langle R^M(T,N)T,B\rangle.
\end{align*}
We have to distinguish the two cases:
\begin{enumerate}
 \item If \(p\neq\frac{1}{2}\) the curve \(\gamma\) is \(p\)-biharmonic if
\begin{align}
\label{system-ode-3d-biharmonic}
k =const\neq 0,\qquad \tau' =\langle R^M(T,N)T,B\rangle,\qquad k ^2+\tau ^2=\langle R^M(T,N)N,T\rangle.
\end{align}
\item If \(p=\frac{1}{2}\) the curve \(\gamma\) is \(p\)-biharmonic if the following
system holds
\begin{align}
\label{system-ode-3d-p12}
0&=\frac{3}{4}k '^2-\frac{1}{2}k k ''-k ^4-\tau ^2k ^2+k ^2\langle R^M(T,N)N,T\rangle, \\
\nonumber0&=-k '\tau +k \tau '-k\langle R^M(T,N)T,B\rangle.
\end{align}
\end{enumerate}
\end{Prop}

Again, we realize that a \(p\)-biharmonic curve in a three-dimensional manifold is biharmonic in general
but for \(p=\frac{1}{2}\) it may also admit additional solutions
corresponding to generalized elastic curves.
Generalized elastic curves in three-dimensional space forms 
have already been classified in
\cite[Section 3]{MR3774309},
but we will again give a short proof 
for the corresponding results on 
proper \(\frac{1}{2}\)-biharmonic curves using our setup, which is
a special case of \cite[Proposition 3.1]{MR3774309}.

\begin{Satz}
\label{thm:curvature-3d}
Let \(\gamma\colon I\to\ M\) be a smooth curve parametrized by arclength and 
\(M\) a three-dimensional Riemannian manifold with constant curvature \(K\).
\begin{enumerate}
\item If  \(M=\s^3\) with the round metric of curvature \(K=1\)
 the curve \(\gamma\) is proper \(\frac{1}{2}\)-biharmonic if
\begin{align}
k(s)=\frac{2}{b_2+\sqrt{b_2^2-4-4b_1^2}\sin 2(s+b_3)},\qquad k(s)=b_1\tau(s),
\end{align}
where \(b_1,b_2,b_3\in\R\).
\item If \(M=\R^3\) with the flat metric the curve \(\gamma\) 
is proper \(\frac{1}{2}\)-biharmonic if
\begin{align}
k(s)=\frac{b_4}{b_4^2(s-b_5)^2+1+b_1^2},\qquad k(s)=b_1\tau(s),
\end{align}
where \(b_1,b_4,b_5\in\R\).
\item If \(M=\mathbb{H}^3\) with the metric of curvature \(K=-1\)
the curve \(\gamma\) is proper \(\frac{1}{2}\)-biharmonic if
\begin{align}
k(s)=\frac{4e^{-2(b_7+s)}}{(b_6+e^{-2(b_7+2)})^2+4+4b_1^2},\qquad k(s)=b_1\tau(s),
\end{align}
where \(b_1,b_6,b_7\in\R\).
\end{enumerate}
\end{Satz}

\begin{Bem}
As in the two-dimensional case the constants \(b_3,b_5,b_7\) can be set to zero
by shifting the curve parameter \(s\).
\end{Bem}

\begin{proof}
By assumption the manifold \(M\) is a three-dimensional space form 
such that
the second equation of \eqref{system-ode-3d-p12} can be written as
\begin{align*}
(\log k )'=(\log\tau )'
\end{align*}
and is solved by
\begin{align*}
k(s)=b_1\tau(s) ,
\end{align*}
where \(b_1\in\R\) is an integration constant.
Then, assuming that we are looking for a non-geodesic solution, 
the first equation of \eqref{system-ode-3d-p12} simplifies to
\begin{align*}
\frac{3}{4}k '^2-\frac{1}{2}k k''-(1+b_1^2)k ^4+Kk^2.
\end{align*}
As in the proof of Theorem \ref{thm:curvature-surface} it is straightforward to see
that the first equation of \eqref{system-ode-3d-p12} is equivalent to
\begin{align*}
\big(k ^{-\frac{1}{2}}\big)''-(1+b_1^2)k^\frac{3}{2}+Kk ^{-\frac{1}{2}}=0.
\end{align*}

Setting \(f:=k ^{-\frac{1}{2}}\) and multiplying by \(f'\) we find 
\begin{align*}
0=(f'^2)'+K(f^2)'+(1+b_1^2)\big(\frac{1}{f^2}\big)',
\end{align*}
which we integrate to obtain
\begin{align*}
b_2=f'^2+Kf^2+\frac{1+b_1^2}{f^2}
\end{align*}
for some \(b_2\in\R\).

To further simplify this equation we set \(h(s):=\frac{1}{2}f^2(s)\)
and get
\begin{align*}
0=h'^2+4Kh^2-2b_2h+1+b_1^2.
\end{align*}
This equation can now easily be integrated distinguishing between the cases 
\(K=1,0,-1\)
and then using that \(k(s)=\frac{1}{2h(s)}\).
\end{proof}

\subsection{Stability of \texorpdfstring{$p$}{}-biharmonic curves on
three-dimensional space forms}
After having determined the structure of \(\frac{1}{2}\)-biharmonic curves
on three-dimensional space forms we now turn to the analysis of their stability
where we again focus on the stability with respect to normal variations.

\begin{Prop}
Let \(\gamma\colon I\to M\) be a \(p\)-biharmonic curve parametrized by arclength
and suppose that \(M\) is a three-dimensional Riemannian space
form of constant curvature \(K\).
Then the Hessian of \eqref{energy-p-biharmonic} evaluated
in the normal direction \(N\) is given by
\begin{align}
\label{eq:second-variation-general-3d}
\operatorname{Hess}(E_{2,p}(\gamma))(N,N)=
\int_I&\big(
(p-1)k^{p-2}(k^2+\tau^2)^2+(k'^2+\tau'^2)k^{p-2}+(p-1)K^2k^{p-2}\\
&-(2p+2)Kk^p
+(2-2p)K\tau^2k^{p-2}
\nonumber\big)ds,
\end{align}
where \(k\) is the geodesic curvature 
and \(\tau\) the torsion
of the curve \(\gamma\).
\end{Prop}

\begin{proof}
We  apply the second variation 
formula derived in Proposition \ref{prop-second-variation-spaceform},
choose \(N\) as variational vector field and 
employ the Frenet-equations \eqref{eq:frenet-3d}.

Inserting this choice into \eqref{eqn:second-variation-space-form}
we find as in the two-dimensional case
\begin{align*}
\operatorname{Hess}(E_{2,p}(\gamma))(N,N) 
=&\nonumber(p-2)\int_I
|\nabla_{\gamma'}\gamma'|^{p-4}
|\langle\nabla_{\gamma'}\nabla_{\gamma'}N,\nabla_{\gamma'}\gamma'\rangle|^2
 ds \\
&\nonumber
 +\int_I|\nabla_{\gamma'}\gamma'|^{p-2}|\nabla_{\gamma'}\nabla_{\gamma'}N|^2ds 
+(p-2)K^2\int_I|\nabla_{\gamma'}\gamma'|^{p-4} 
|\langle\nabla_{\gamma'}\gamma',N\rangle|^2 ds \\ 
&\nonumber+K^2\int_I|\nabla_{\gamma'}\gamma'|^{p-2}ds
\nonumber+K\int_I
\langle\gamma',\nabla_{\gamma'}\big(|\nabla_{\gamma'}\gamma'|^{p-2}\nabla_{\gamma'}\gamma'\big)\rangle
ds \\
&\nonumber
+3K\int_I|\nabla_{\gamma'}\gamma'|^{p-2}
\langle\nabla_{\gamma'}N,\gamma'\rangle
\langle\nabla_{\gamma'}\gamma',N\rangle ds  \\
&\nonumber
+2(p-2)K\int_I
\langle\nabla_{\gamma'}\nabla_{\gamma'}N,\nabla_{\gamma'}\gamma'
\rangle 
|\nabla_{\gamma'}\gamma'|^{p-4}
\langle N,\nabla_{\gamma'}\gamma'\rangle
 ds\\
&\nonumber+2K\int_I|\nabla_{\gamma'}\gamma'|^{p-2}
\langle\nabla_{\gamma'}\nabla_{\gamma'}N,N\rangle
ds.
\end{align*}

Using the Frenet-equations \eqref{eq:frenet-3d} it is straightforward
to calculate that
\begin{align*}
\langle\nabla_{\gamma'}\nabla_{\gamma'}N,\nabla_{\gamma'}\gamma'\rangle
=&-k(k^2+\tau^2),\\
|\nabla_{\gamma'}\nabla_{\gamma'}N|^2
=&k'^2+\tau'^2+(k^2+\tau^2)^2.
\end{align*}
In addition, we find
\begin{align*}
\int_I
\langle\gamma',\nabla_{\gamma'}\big(|\nabla_{\gamma'}\gamma'|^{p-2}\nabla_{\gamma'}\gamma'\big)\rangle
ds =&
\int_I|\nabla_{\gamma'}\gamma'|^{p-2}
\langle\nabla_{\gamma'}N,\gamma'\rangle
\langle\nabla_{\gamma'}\gamma',N\rangle ds  
=-\int_Ik^pds, \\
\int_I
\langle\nabla_{\gamma'}\nabla_{\gamma'}N,\nabla_{\gamma'}\gamma'
\rangle 
|\nabla_{\gamma'}\gamma'|^{p-4}
\langle N,\nabla_{\gamma'}\gamma'\rangle
 ds=&-\int_I(k^2+\tau^2)k^{p-2}ds,\\
\int_I|\nabla_{\gamma'}\gamma'|^{p-2}
\langle\nabla_{\gamma'}\nabla_{\gamma'}N,N\rangle
ds=&-\int_I(k^2+\tau^2)k^{p-2}ds.
\end{align*}
Inserting into the above equation completes the proof.
\end{proof}

\begin{Satz}
Let \(\gamma\colon I\to M\) be a \(p\)-biharmonic curve 
parametrized by arclength
and suppose that \(M\) is a three-dimensional 
Riemannian manifold of constant curvature \(K\).
Consider the second variation of \eqref{energy-p-biharmonic}
for a variation in the normal direction \(N\). Then the following
statements hold:
\begin{enumerate}
\item For \(p\neq\frac{1}{2}\)
  we get 
  \begin{align}
  \label{hessian-p-general-3d}
\operatorname{Hess}(E_{2,p}(\gamma))(N,N) 
=-4K|I|k^{p}.
 \end{align}
\item For \(p=\frac{1}{2}\) we get
\begin{align}
\label{hessian-p12-3d}
\operatorname{Hess}&(E_{2,p}(\gamma))(N,N) \\
\nonumber &=\int_I\big(
-\frac{(1+a^2)^2}{2}k^\frac{5}{2}+(1+a^2)k'^2k^{-\frac{3}{2}}
-\frac{1}{2}K^2k^{-\frac{3}{2}}-3Kk^\frac{1}{2}
+a^2Kk^\frac{1}{2}
\big)ds,
\end{align}
where \(a\in\R\).

 \end{enumerate}
\end{Satz}
\begin{proof}
First, we consider the case \(p\neq\frac{1}{2}\).
Then, from \eqref{system-ode-3d-biharmonic} we know that
the condition for a curve to be \(p\)-biharmonic 
is \(k^2+\tau^2=K\). Inserting into \eqref{eq:second-variation-general-3d}
then completes the proof for this case.

Concerning the case \(p=\frac{1}{2}\)
we know that the geodesic curvature and the torsion of the curve
are related by 
\(\tau(s)=ak(s)\)
with \(a\in\R\).
Inserting into \eqref{eq:second-variation-general-3d} we get the result.
\end{proof}

\begin{Bem}
\begin{enumerate}
 \item Note that for \(p\neq\frac{1}{2}\) the torsion \(\tau\) of the \(p\)-biharmonic curve enters into \eqref{hessian-p-general-3d} via \(K=k^2+\tau^2\).
 Consequently, we get the same formula as in the two-dimensional case, that is \eqref{hessian-p-general} and \eqref{hessian-p-general-3d} coincide.
 \item The same statement as made in Remark \ref{rem-yelin} also holds
 in the case of a three-dimensional target \(M\).
  \item As in the surface case  we can insert the solutions
  obtained in Theorem \ref{thm:curvature-3d} for \(p=\frac{1}{2}\) into equation
 \eqref{hessian-p12-3d} in order to check the sign of the second variation.
 Unfortunately, this again leads to a number of complicated expressions
 such that no information can be extracted.
\end{enumerate}

\end{Bem}

\par\medskip
\textbf{Acknowledgements:}
The author would like to thank the referee both for carefully reading the manuscript and for many valuable suggestions.
\bibliographystyle{plain}
\bibliography{mybib}

\begin{thebibliography}{10}

\bibitem{MR2031956}
J.~Arroyo, O.~J. Garay, and M.~Barros.
\newblock Closed free hyperelastic curves in the hyperbolic plane and
  {C}hen-{W}illmore rotational hypersurfaces.
\newblock {\em Israel J. Math.}, 138:171--187, 2003.

\bibitem{MR2007599}
J.~Arroyo, O.~J. Garay, and J.~J. Menc\'{\i}a.
\newblock Closed generalized elastic curves in {$S^2(1)$}.
\newblock {\em J. Geom. Phys.}, 48(2-3):339--353, 2003.

\bibitem{MR2078687}
J.~Arroyo, O.~J. Garay, and J.~J. Menc\'{\i}a.
\newblock Extremals of curvature energy actions on spherical closed curves.
\newblock {\em J. Geom. Phys.}, 51(1):101--125, 2004.

\bibitem{MR3774309}
J.~Arroyo, O.~J. Garay, and A.~P\'{a}mpano.
\newblock Constant mean curvature invariant surfaces and extremals of curvature
  energies.
\newblock {\em J. Math. Anal. Appl.}, 462(2):1644--1668, 2018.

\bibitem{MR4028018}
J.~Arroyo, O.~J. Garay, and A.~P\'{a}mpano.
\newblock Delaunay surfaces in {$\Bbb S^3(\rho)$}.
\newblock {\em Filomat}, 33(4):1191--1200, 2019.

\bibitem{MR1548296}
Wilhelm Blaschke.
\newblock Literaturberichte: {V}orlesungen \"{u}ber {D}ifferentialgeometrie.
\newblock {\em Monatsh. Math. Phys.}, 23(1):A12--A13, 1912.
\newblock Von L. Bianchi. Autorisierte deutsche \"{U}bersetzung von M. Lukat.
  Zweite, vermehrte und verbesserte Auflage. Leipzig, B. G. Teubner 1910.

\bibitem{MR4418802}
Simon Blatt, Christopher Hopper, and Nicole Vorderobermeier.
\newblock A regularized gradient flow for the {$p$}-elastic energy.
\newblock {\em Adv. Nonlinear Anal.}, 11(1):1383--1411, 2022.

\bibitem{MR4106647}
V.~Branding, S.~Montaldo, C.~Oniciuc, and A.~Ratto.
\newblock Higher order energy functionals.
\newblock {\em Adv. Math.}, 370:107236, 60, 2020.

\bibitem{MR4058514}
Volker Branding.
\newblock On interpolating sesqui-harmonic maps between {R}iemannian manifolds.
\newblock {\em J. Geom. Anal.}, 30(1):248--273, 2020.

\bibitem{MR3778117}
Volker Branding and Florian Hanisch.
\newblock Magnetic geodesics via the heat flow.
\newblock {\em Asian J. Math.}, 21(6):995--1014, 2017.

\bibitem{MR4477489}
Volker Branding and Anna Siffert.
\newblock On the equivariant stability of harmonic self-maps of cohomogeneity
  one manifolds.
\newblock {\em J. Math. Anal. Appl.}, 517(2):Paper No. 126635, 19, 2023.

\bibitem{MR2250208}
R.~Caddeo, S.~Montaldo, C.~Oniciuc, and P.~Piu.
\newblock The classification of biharmonic curves of {C}artan-{V}ranceanu
  3-dimensional spaces.
\newblock In {\em Modern trends in geometry and topology}, pages 121--131. Cluj
  Univ. Press, Cluj-Napoca, 2006.

\bibitem{MR1884940}
R.~Caddeo, S.~Montaldo, and P.~Piu.
\newblock Biharmonic curves on a surface.
\newblock {\em Rend. Mat. Appl. (7)}, 21(1-4):143--157, 2001.

\bibitem{MR3567234}
Xiangzhi Cao and Yong Luo.
\newblock On {$p$}-biharmonic submanifolds in nonpositively curved manifolds.
\newblock {\em Kodai Math. J.}, 39(3):567--578, 2016.

\bibitem{MR3098705}
Yuan-Jen Chiang.
\newblock {\em Developments of harmonic maps, wave maps and {Y}ang-{M}ills
  fields into biharmonic maps, biwave maps and bi-{Y}ang-{M}ills fields}.
\newblock Frontiers in Mathematics. Birkh\"{a}user/Springer, Basel, 2013.

\bibitem{MR3781340}
Anna Dall'Acqua, Tim Laux, Chun-Chi Lin, Paola Pozzi, and Adrian Spener.
\newblock The elastic flow of curves on the sphere.
\newblock {\em Geom. Flows}, 3(1):1--13, 2018.

\bibitem{MR4663707}
Anthony Gruber, \'{A}lvaro P\'{a}mpano, and Magdalena Toda.
\newblock Instability of closed {$p$}-elastic curves in {$(\Bbb{S})^2$}.
\newblock {\em Anal. Appl. (Singap.)}, 21(6):1533--1559, 2023.

\bibitem{MR4169439}
Yingbo Han and Yong Luo.
\newblock Nonexistence of proper {$p$}-biharmonic maps and {L}iouville type
  theorems {I}: case of {$p\ge2$}.
\newblock {\em J. Elliptic Parabol. Equ.}, 6(2):409--426, 2020.

\bibitem{MR3393120}
Yingbo Han and Wei Zhang.
\newblock Some results of {$p$}-biharmonic maps into a non-positively curved
  manifold.
\newblock {\em J. Korean Math. Soc.}, 52(5):1097--1108, 2015.

\bibitem{MR2076749}
Rongpei Huang.
\newblock A note on the {$p$}-elastica in a constant sectional curvature
  manifold.
\newblock {\em J. Geom. Phys.}, 49(3-4):343--349, 2004.

\bibitem{MR886529}
Guo~Ying Jiang.
\newblock {$2$}-harmonic maps and their first and second variational formulas.
\newblock {\em Chinese Ann. Math. Ser. A}, 7(4):389--402, 1986.
\newblock An English summary appears in Chinese Ann. Math. Ser. B {{\bf{7}}}
  (1986), no. 4, 523.

\bibitem{MR3773492}
S.~M. Kazemi~Torbaghan and M.~M. Rezaii.
\newblock Stability of {F}-biharmonic maps.
\newblock {\em Bull. Iranian Math. Soc.}, 43(6):1657--1669, 2017.

\bibitem{MR772124}
Joel Langer and David~A. Singer.
\newblock The total squared curvature of closed curves.
\newblock {\em J. Differential Geom.}, 20(1):1--22, 1984.

\bibitem{MR4277362}
Carlo Mantegazza, Alessandra Pluda, and Marco Pozzetta.
\newblock A survey of the elastic flow of curves and networks.
\newblock {\em Milan J. Math.}, 89(1):59--121, 2021.

\bibitem{MR4594934}
Ahmed Mohammed~Cherif and Khadidja Mouffoki.
\newblock {$p$}-biharmonic hypersurfaces in {E}instein space and conformally
  flat space.
\newblock {\em Bull. Korean Math. Soc.}, 60(3):705--715, 2023.

\bibitem{MR4216418}
S.~Montaldo, C.~Oniciuc, and A.~Ratto.
\newblock Index and nullity of proper biharmonic maps in spheres.
\newblock {\em Commun. Contemp. Math.}, 23(3):Paper No. 1950087, 36, 2021.

\bibitem{MR3045700}
Stefano Montaldo and Andrea Ratto.
\newblock A general approach to equivariant biharmonic maps.
\newblock {\em Mediterr. J. Math.}, 10(2):1127--1139, 2013.

\bibitem{MR4150936}
Shinya Okabe, Paola Pozzi, and Glen Wheeler.
\newblock A gradient flow for the {$p$}-elastic energy defined on closed planar
  curves.
\newblock {\em Math. Ann.}, 378(1-2):777--828, 2020.

\bibitem{MR4110268}
Ye-Lin Ou.
\newblock A note on equivariant biharmonic maps and stable biharmonic maps.
\newblock {\em J. Math. Anal. Appl.}, 491(1):124301, 11, 2020.

\bibitem{MR4386842}
Ye-Lin Ou.
\newblock Stability and the index of biharmonic hypersurfaces in a {R}iemannian
  manifold.
\newblock {\em Ann. Mat. Pura Appl. (4)}, 201(2):733--742, 2022.

\bibitem{MR4265170}
Ye-Lin Ou and Bang-Yen Chen.
\newblock {\em Biharmonic submanifolds and biharmonic maps in {R}iemannian
  geometry}.
\newblock World Scientific Publishing Co. Pte. Ltd., Hackensack, NJ, [2020]
  \copyright 2020.

\bibitem{MR4318851}
Marco Pozzetta.
\newblock Convergence of elastic flows of curves into manifolds.
\newblock {\em Nonlinear Anal.}, 214:Paper No. 112581, 53, 2022.

\bibitem{MR2788659}
Matthias Schneider.
\newblock Closed magnetic geodesics on {$S^2$}.
\newblock {\em J. Differential Geom.}, 87(2):343--388, 2011.

\bibitem{MR2959932}
Matthias Schneider.
\newblock Closed magnetic geodesics on closed hyperbolic {R}iemann surfaces.
\newblock {\em Proc. Lond. Math. Soc. (3)}, 105(2):424--446, 2012.

\bibitem{MR4184824}
Naoki Shioji and Kohtaro Watanabe.
\newblock Total {$p$}-powered curvature of closed curves and flat-core closed
  {$p$}-curves in {${\rm S}^2(G)$}.
\newblock {\em Comm. Anal. Geom.}, 28(6):1451--1487, 2020.

\bibitem{MR3676571}
Masaaki Umehara and Kotaro Yamada.
\newblock {\em Differential geometry of curves and surfaces}.
\newblock World Scientific Publishing Co. Pte. Ltd., Hackensack, NJ, 2017.
\newblock Translated from the second (2015) Japanese edition by Wayne Rossman.

\bibitem{MR3229087}
Kohtaro Watanabe.
\newblock Planar {$p$}-elastic curves and related generalized complete elliptic
  integrals.
\newblock {\em Kodai Math. J.}, 37(2):453--474, 2014.

\end{thebibliography}
\end{document}